\documentclass[a4paper, twoside, 12pt]{amsart}

\usepackage{graphicx}
\usepackage{framed}

\usepackage[utf8]{inputenc}

\usepackage[all]{xy}
\usepackage{pb-diagram, pb-xy}

\usepackage{amssymb}
\usepackage{amsthm}

\usepackage{mathptmx}

\numberwithin{equation}{section}

\theoremstyle{plain}
\newtheorem*{cor*}{Corollary}
\newtheorem*{ass*}{Assumption}

\newcommand{\bbC}{{\mathbb C}}

\newcommand{\bbQ}{{\mathbb Q}}

\newcommand{\calH}{{\mathcal H}}

\newcommand{\bfD}{{\mathbf{D}}}

\begin{document}

\title[A Remark on Rigidity of BN-sheaves]{A Remark on Rigidity of BN-sheaves}
\author{R. Weissauer}

\maketitle

\thispagestyle{empty}


Let $X$ be an abelian variety over
an algebraically closed field $k$, assuming that either $k=\bbC$ or that $k$ is the algebraic closure of a finite field.  In [BN] we considered the convolution
product $K*L$ for complexes $K$ and $L$ in the bounded 
derived category $\bfD=D_c^b(X,\Lambda)$ with coefficient field $\Lambda$
either $\bbC$
for $k=\bbC$, or $\Lambda=\overline\bbQ_l$. The convolution product 
is defined, using the group law $a: X\times X\to X$ of the abelian 
variety $X$, as the derived direct image complex $K*L=Ra_*(K\boxtimes L)$. 
This convolution product makes $(\bfD,*)$ into a symmetric monoidal category; 
for details we refer
to [BN]. For a complex $K$ let $D(K)$ denote its Verdier dual.

\medskip
 In a symmetric
monoidal category $(\bfD,\otimes)$ an object $M$ is called rigid, if the 
functor $F(K)=K\otimes M$ admits as right adjoint functor $G(N)=M^\vee \otimes N$ for some object $M^\vee$ in $\bfD$, which then is
called the dual of $M$.
We remark, that in the literature rigidity is often stated differently
in terms of evaluation and coevaluation morphisms. However this 
alternative definition is equivalent to the one given above
as an easy consequence of
[ML, IV.1, thm. 1]. We now
generalize [BN, cor.6] with a rather simple argument.

\bigskip\noindent
{\bf Proposition}. {\it Any $M\in D_c^b(X,\Lambda)$ is rigid 
with dual $M^\vee = (-id_X)^*D(M)$.
  }
  
\bigskip\noindent
{\it Proof}. Adjoint functors $F:\bfD\to \bfD'$ and $G:\bfD'\to \bfD$ between 
categories $\bfD$ and $\bfD'$ are given by the following data: A collection of
isomorphisms 
$$ \varphi_{K,N}: Hom_{\bfD}(F(K),N) \cong Hom_{\bfD'}(K,G(N)) \ ,$$ 
for each $K\in \bfD$, $N\in \bfD'$, which are functorial in $K$ and $N$. 
Let $\bfD=\bfD'$ be $D_c^b(X,\overline\bbQ_l)$. For our proof it is enough to
show that $F(K)=K*M$ admits $G(N)=M^\vee*N$ as an adjoint for $M^\vee =  (-id_X)^*D(M)$. For this
we have to construct functorial isomorphisms $\varphi_{K,N}$ as above.

\bigskip\noindent
Let $b:X\to Spec(k)=\{0\}$ be the structure morphism, let $a$ be the group law and $f:X\times X\to X$ be the difference morphism $f(x,y)=y-x$. All these morphisms are proper, hence $Ra_*=Ra_!, Rb_*=Rb_!$ etc. By unravelling the definitions for $L_0 := Rb_*(Ra_*(K\boxtimes M)\otimes^L D(N))$ we get
$$ Hom_{\bfD}(K*M,N)=Hom_{\bfD}(Ra_*(K\boxtimes M),N) $$ $$ =Rb_*(R{\calH}om(Ra_*(K\boxtimes M),N)) $$
$$ = Rb_*(D(Ra_*(K\boxtimes M)\otimes^L D(N))) $$ $$ = 
D(Rb_!(Ra_*(K\boxtimes M)\otimes^L D(N)))
 = D(L_0) \ .$$
Similarly
for $R_0:=  Rb_*(K \otimes^L Rf_*(M \boxtimes D(N))
$ we get $$ Hom_{\bfD}(K,M^\vee*N)=Hom_{\bfD}(K,Ra_*(M^\vee\boxtimes N)) $$ $$
=Rb_*(R{\calH}om(K,Ra_*(M^\vee \boxtimes N))) $$
$$=Rb_*( D(K \otimes^L D(Ra_*(M^\vee \boxtimes N)))) $$ $$
= D(Rb_!(K \otimes^L Ra_!(D(M^\vee \boxtimes N))))) $$
$$ = D(Rb_*(K \otimes^L Ra_*(D(M^\vee) \boxtimes D(N)))) $$ $$
 = D(Rb_*(K \otimes^L Ra_*((-id_X)^*M \boxtimes D(N)) ))
 = D(R_0) \ .$$
Here, by abuse of notation, the equality signs indicate functorial isomorphisms,
whose precise nature however will not be important for us. Indeed, for our proof, it is now enough
to find  isomorphisms 
$$  \psi_{A,C}: \ \ L_0=Rb_*(Ra_*(A\boxtimes B) \otimes^L C) \ \ \cong \ \ R_0 =
Rb_*(A \otimes^L Rf_*(B\boxtimes C)) \ $$
for fixed $B$, functorial in $A$ and $C$. 
For this we use the proper base change theorem. Notice $z-(x+y)=(z-y)-x$ or $f\circ (a\times id_X) = f \circ (id_X \times f)$ holds,
considered as morphisms $X\times X \times X \to X$. Therefore there exist isomorphisms
$$ \Psi_{A,C}: \ L \ \cong \ R $$
functorial  in $A$ and $C$ (for fixed $B$) between
$$ L:=Rf_*(R(a\times id_X)_* (A\boxtimes B \boxtimes C)) 
$$ and $$
R := Rf_*(R(id_X\times f)_*  (A\boxtimes B \boxtimes C)) \ .$$
Hence the stalks $R_{\{0\}}$ resp. $L_{\{0\}}$ of the complexes $R$ resp. $L$ at the origin $0$ can be 
identified. Now via the proper base change theorem, applied for the cartesian diagram
$$ \xymatrix{  X \ar[d]_b \ar[r]^-\Delta & X\times X \ar[d]^f \cr \{0\} \ar[r]^-{i_0} & X \cr } $$
obtained by the diagonal embedding $\Delta(x)=(x,x)$, we can compute  for any complex $K$ in $D_c^b(X\times X,\Lambda)$ the stalk $Rf_*(K)_{\{0\}} \cong Rb_*(\Delta^*(K))$ of the direct image $Rf_*(K)$ at the origin. This implies $R_{\{0\}}  =R_0$, since
$$ R_{\{0\}}  =  Rf_*(R(id_X\times f)_*  (A\boxtimes B \boxtimes C))_{\{0\}}
 =  Rf_*(  A\boxtimes Rf_*(B \boxtimes C))_{\{0\}} $$
$$ = Rb_*( \Delta^*(A \boxtimes Rf_*(B\boxtimes C)))  =
Rb_*(A\otimes^L Rf_*(B\boxtimes C)) = R_0 \ ,$$
and similarly $ L_{\{0\}} = L_0$. Thus the stalk morphisms $\psi_{A,C} = (\Psi_{A,C})_{\{0\}}$ give the desired functorial isomorphisms we have to construct. \qed
 
\medskip 
{\bf References}

\medskip [ML] MacLane S., {\em Categories for the working mathematician}, Springer  

\medskip [BN] Weissauer R., {\em Brill-Noether sheaves}, arXiv:math/0610923v4
 
\end{document}